\documentclass[a4paper]{article}

\usepackage{amsmath,amssymb,amsthm} 
\usepackage[english]{babel}

\newtheorem{thm}{Theorem}[section]

\newtheorem{lemma}[thm]{Lemma}

\theoremstyle{definition}
\newtheorem{definition}[thm]{Definition}
\newtheorem{example}[thm]{Example}
\theoremstyle{remark}

\numberwithin{equation}{subsection}

\def\demo{\noindent \textit{Proof: }}
\def\Biy{\mathrm{Biy}}

\def\Mot{\mathrm{Mot}}
\def\AA{\mathbb{A}}
\def\PP{\mathbb{P}}
\def\PGl{\mathrm{PGl}}
\def\RR{\mathbb{R}}

\begin{document}

\title{\textsc{Automorphisms of classical geometries in the sense of Klein}}

\author{Navarro, A.\thanks{ICMAT, Madrid, Spain}
\and Navarro, J.\thanks{Corresponding author. Department of
Mathematics, University of Extremadura, Badajoz, Spain.\newline {\it
Email address:} navarrogarmendia@unex.es}}


\maketitle

\begin{abstract}
In this note, we compute the group of  automorphisms of Projective,
Affine and Euclidean Geometries in the sense of Klein.

As an application, we give a simple construction of the outer automorphism of $S_6$.
\end{abstract}

\bigskip


\section{Introduction}

Let $\mathbb{P}_n$ be the set of 1-dimensional subspaces of a
$n+1$-dimensional vector space $E$ over a (commutative) field $k$.
This standard definition does not capture the "structure" of the
projective space although it does point out its automorphisms: they
are projectivizations of semilinear automorphisms of $E$, also known
as {\it Staudt projectivities}.

A different approach (v. gr. \cite{Artin}), defines the projective
space as a lattice (the lattice of all linear subvarieties)
satisfying certain axioms. Then, the so named Fundamental Theorem of
Projective Geometry (\cite{Artin}, Thm 2.26) states that, when $n >
1$, collineations of $\mathbb{P}_n$ (i.e., bijections preserving
alignment, which are the automorphisms of this lattice structure)
are precisely Staudt projectivities.

In this note we are concerned with geometries in the sense of Klein:
a \textit{geometry} is  a pair $(X, G)$ where $X$ is a set and $G$
is a subgroup of the group $\Biy (X)$ of all bijections of $X$. In
Klein's view, Projective Geometry is the pair $(\PP_n ,
\mathrm{PGl}_n$), where $\mathrm{PGl}_n$ is the group of
projectivizations of $k$-linear automorphisms of the vector space
$E$ (see Example \ref{ExProj}). The main result of this note is a
computation, analogous to the aforementioned theorem for
collineations, but in the realm of Klein geometries:

\medskip
\noindent {\bf Theorem \ref{fundamental}} {\it The group of
automorphisms of the Projective Geometry
$(\mathbb{P}_n , \mathrm{PGl}_n)$ is the group of
Staudt projectivities, for any $n \geq 1$.}
\medskip

Let us remark this statement covers the case $n=1$ of the projective
line, which usually requires a separated treatment (\cite{Artin},
\cite{Sancho}).

Also, we compute in Theorems \ref{afin} and \ref{Euclidean} the
group of automorphisms for other classical geometries -- namely
Affine and Euclidean Geometries.

Finally, as an application, we use Theorem \ref{fundamental} to give
a description of the outer automorphism of $S_6$, which is simple in
comparison with other constructions existing in the literature (see
\cite{Invariants, Lotman} or \cite{Lam}).

\section{Preliminaries}

Let us firstly introduce the main definitions and examples that will be used later on:

\begin{definition}
A \textit{geometry} in the sense of Klein is a pair $(X,G)$ where
$X$ is a set and $G$ is a subgroup of the group $\Biy (X)$ of all
bijections $X\to X$.

\end{definition}

The {\it concepts} of a geometry are those notions invariant with
respect to the action of the structural group $G$.

\begin{example}\label{ExProj}
The examples that will appear in this note are:

\begin{itemize}
\item  The \emph{Projective Geometry} $(\mathbb{P}_n , \mathrm{PGl}_n)$,
where $\mathbb{P}_n$ is the set of 1--dimensional vector subspaces
of a $k$--vector space $E$ of dimension $n+1$ and $\mathrm{PGl}_n$
is the group projectivizations of $k$--linear automorphisms of $E$
with the obvious action on $\mathbb{P}_n$.

\item The \emph{Affine Geometry} $(\mathbb{A}_n,
\mathrm{Aff}_n)$, where $\mathbb{A}_n$ is the complement of an
hyperplane $H$ in $\mathbb{P}_n$, called \emph{hyperplane at
infinity}, and $\mbox{Aff}_n$ is the subgroup of $\mbox{PGl}_n$
consisting of all projectivities $\varphi \colon \PP_n \to \PP_n$
preserving the hyperplane at infinity (i.e., such that $\varphi (H)
= H$), with the action by restriction on $\AA _n$.

\item The \emph{Euclidean Geometry} $(\AA_n(\mathbb{R})$, $\mbox{Mot}_n)$,
where $\AA_n (\RR )$ denotes the affine space over the real numbers,
endowed with a positive definite, non-singular metric $g$ on its
vector space of directions, and $\mbox{Mot}_n$ stands for the group
of {\it motions}; i.e., the group of those affinities $\varphi
\colon \AA_n (\RR) \to \AA_n (\RR)$ whose tangent linear map
$\vec{\varphi}$ satisfy $\vec{\varphi}_* g = g$.

\end{itemize}
\end{example}

\begin{definition}

An {\it isomorphism} between two geometries $(X,G)$ and $(X',G')$ is
a bijection $f \colon X\to X'$ such that the map $$\phi \colon
\mbox{Biy}(X) \to \mbox{Biy}(X') \qquad , \qquad  \phi (g) :=
\varphi \circ g \circ \varphi ^{-1}$$  preserves the structural
groups, $\phi (G) = G' $; that is to say, such that $\phi G \phi
^{-1}=G'$.

An {\it automorphism} of a geometry $(X, G)$ is an isomorphism of
geometries $f \colon (X,G)\to (X,G)$. With the composition of maps,
the automorphisms of a geometry $(X, G)$ are a group.
\end{definition}

In other words, the group of automorphisms of a geometry $(X,G)$ is
the normalizer of $G$ inside the group $\Biy (X)$ of all bijections
of $X$.

\section{Automorphisms of classical geometries}

\subsection{Projective Geometry}

Let $\mathbb{P}_n$ be the set of 1-dimensional subspaces of a
$(n+1)$-dimensional $k$-vector space $E$ and let:
$$
\pi \colon E {\rm -} \{ 0 \} \longrightarrow \PP_n \quad , \quad e
\longmapsto \langle e \rangle \ ,
$$
be the  projectivization map. Recall that linear subvarieties of
$\PP_n $ (lines, planes,...) are defined as the projectivization of
linear subspaces of $E$.

\begin{definition}
Three points $p_1, p_2, p_3 \in \PP_n$ are {\it collinear} if there
exists a line $L$ passing through them.

A bijection $\varphi \colon \PP_n \to \PP_n$ is a {\it collineation}
if it transforms lines into lines; that is, it maps collinear points
into collinear points.
\end{definition}

\begin{definition}
A {\it semilinear automorphism} $f\colon E\to E$ is a bijection such
that there exists an automorphism of fields $h\colon k\to k$
satisfying
$$f(\lambda e + \mu v)=
h(\lambda) f(e)+h(\mu) f(v) \quad , \quad \forall \ e,v \in E, \  \
\lambda, \mu \in k \ . $$

A bijection $\Phi \colon \mathbb{P}_n\to \mathbb{P}_n$ is a {\it
Staudt projectivity} if it is the projectivization of a semilinear
automorphism.
\end{definition}

Observe that elements in $\, \PGl_n \,$ are precisely the Staudt
projectivities with associated automorphism of fields $h = {\rm
Id}$.

\medskip

For what follows, it will be useful to characterize
colli\-nea\-tions in terms of projectivities; to this end, let us
consider the following sets:
$$ \mathsf{P}_{p_1,p_2}(p_3) := \{ \,
p \in \PP _n {\rm - } \{ p_3 \} \colon \ \exists \ \varphi \in
\mathrm{PGl}_n \ , \, \varphi (p_1)=p_1,\ \varphi (p_2)=p_2,\
\varphi (p_3)= p \, \} \ . $$

\begin{lemma}\label{alineacion}
Three different points $p_1,p_2,p_3 \in \mathbb{P} _n$ are collinear
if and only if
$$
\mathsf{P}_{p_1,p_2}(p_3)= \mathsf{P}_{p_1,p_3}(p_2) =
\mathsf{P}_{p_2,p_3}(p_1) \ . $$ \end{lemma}

\demo If $p_1,p_2,p_3$ are collinear, then
$\mathsf{P}_{p_1,p_2}(p_3)$ equals the complement of $p_1,p_2,p_3$
in the line passing through them.

Conversely, if $p_1,p_2,p_3$ are not collinear, the set
$\mathsf{P}_{p_1,p_2}(p_3)$ is the complement of the line $L$
joining $p_1$ and $p_2$ (apart from $p_3$). Then, any point $q\in
L$, $q\neq p_1,p_2$ satisfies that $q\notin
\mathsf{P}_{p_1,p_2}(p_3)$ and $q \in \mathsf{P}_{p_1,p_3}(p_2) \cap
\mathsf{P}_{p_2,p_3}(p_1)$.

\hfill $\square$

\begin{thm}\label{fundamental}
The group of automorphisms of the Projective Geometry
$(\mathbb{P}_n,\mathrm{PGl}_n)$ is the group of
Staudt projectivities, for any $n \geq 1$.

That is to say, the group of Staudt projectivities is the normalizer
of the group $\mathrm{PGl}_n$ in the group of all bijections of
$\mathbb{P}_n$.
\end{thm}

\demo Let $\Phi \colon \PP_n \to \PP_n\,$ be a Staudt projectivity,
with associated automorphism of fields $h \colon k \to k$.

If $\varphi \in \PGl_n$ then $\Phi \varphi \Phi^{-1}$ is a Staudt
projectivity, with associated automorphism $h\circ \mathrm{Id}\circ
h^{-1}=\mathrm{Id}$, so that $\Phi \mathrm{PGl}_n \Phi^{-1}\subseteq
\mathrm{PGl}_n$. Since $\Phi^{-1}$ is also a Staudt projectivity, it
also holds $\Phi^{-1} \mathrm{PGl}_n \Phi \subseteq \mathrm{PGl}_n$
and the reverse inclusion follows. We conclude that $\Phi$ belongs
to the normalizer of $\mathrm{PGl}_n$ in the group
$\mathrm{Biy}(\mathbb{P}_n)$.

\medskip

Now, let us prove that any automorphism $\varphi$ of
$(\mathbb{P} _n , \mathrm{PGl}_n)$ is a Staudt projectivity.

\medskip

$n\geq 2.$ As we already mentioned, in this case Staudt
projectivities are precisely collineations (\cite{Artin}, Thm 2.26),
and hence it is enough to check that $\varphi$ is a collineation.
Since $\varphi$ preserves elements in $\,\PGl_n\,$ by hypothesis, it
easily follows that
$$\varphi (\mathsf{P}_{p_1,p_2}(p_3))=
\mathsf{P}_{\varphi (p_1), \varphi (p_2)} (\varphi (p_3)) \ , $$ so
Lemma \ref{alineacion} allows to deduce that $\varphi (p_1), \varphi
(p_2),\varphi (p_3)$ are collinear whenever so are $p_1,p_2,p_3$.

\medskip

$n=1.$ Let us fix a projective reference $(p_0, p_\infty, p_1)$ in
$\mathbb{P}_1$, and write $ p'_0 := \varphi (p_0), \ p'_1 :=\varphi
(p_1), \ p'_\infty := \varphi (p_\infty)$.

In the affine line $\mathbb{A}_1 = \mathbb{P}_1 {\rm -}\{ p_\infty
\} $, consider the origin $p_0$ and the unit point $p_1$, thus
inducing a bijection $\mathbb{A}_1 \simeq k$. Let $p_\lambda$ denote
the point corresponding to $\lambda \in k$ (and analogously in the
affine line $\mathbb{A}'_1 = \mathbb{P}_1 {\rm -}\{ p'_\infty \} $).

The composition $k \simeq \mathbb{A}_1 \xrightarrow{\varphi}
\mathbb{A}'_1 \simeq k $ defines a bijection $h \colon k\to k$ that
is an automorphism of the field $k$ (Lemma \ref{automorfismo}) and
such that $\varphi(p_\lambda)=p'_{ h(\lambda)}$.

Therefore, if $(e_0,e_1)$ is a basis of $E$ normalized to the
reference $(p_0,p_\infty , p_1)$ and $(e'_0,e'_1)$ is a basis
normalized to the reference $(p'_0 , p'_\infty , p'_1)$, then
$\varphi$ coincides with the projectivization of the semilinear map
$f( \lambda e_0+ \mu e_1)=h(\lambda)\, e'_0+h(\mu)\, e'_1$.

\hfill
$\square$

\begin{lemma}\label{automorfismo}
The map $h\colon k \rightarrow k$ defined in the proof above is an
automorphism of fields.
\end{lemma}

\demo By definition, $h(0)=0$ and $h(1)=1$, so we only have to show
that $h$ is compatible with additions and products.

Let $\tau_{\mu}$ and $\sigma_{\mu}$ denote the translation by $\mu
\in k$ in $\mathbb{A}_1$ and the homothety with center at $p_0$ and
ratio $\mu$, respectively (and analogously, $\tau_{\mu}'$ and
$\sigma_{\mu}'$ in $\mathbb{A}_1'$). \medskip

Let us first check the equality $\varphi \tau_{\mu} \varphi ^{-1} =
\tau'_{h(\mu)}$: since $p_\infty$ is the only point that is fixed by
$\tau_{\mu}$, its image $p'_\infty$ is the unique fixed point of
$\varphi \tau_{\mu} \varphi ^{-1}$; hence this composition is a
translation in $\mathbb{A}_1'$ (it is an homography by hypothesis),
and it suffices to see that it transforms $p'_0$ into $p'_{h(\mu)}$:

$$ \varphi \tau_{\mu} \varphi ^{-1} (p'_0)= \varphi
\tau_{\mu}(p_0)= \varphi (p_\mu)= p'_{h(\mu)}.$$

Using this equality, $\varphi \tau_{\mu} = \tau'_{h(\mu)}\varphi$,
and applying it to $p_\lambda$:
\begin{align*}
& \varphi \tau_{\mu} (p_\lambda)=\varphi(p_{\lambda + \mu})= p'_{h(\lambda + \mu)} \\
& \tau'_{h(\mu)} \varphi
(p_\lambda)=\tau'_{h(\mu)}(p'_{h(\lambda)}) = p'_{h(\lambda) +
h(\mu)} \end{align*} it follows that $h(\lambda + \mu)=
h(\lambda) + h(\mu)$.

\medskip
In a similar way, $ \varphi \sigma_{\mu} \varphi ^{-1}$ can be
proved to be the homothety $\sigma ' _{h(\mu)}$: since $p_0$ and
$p_\infty$ are fixed points of $\sigma_\mu$, both $p'_0$ and
$p'_\infty$ are fixed points of $ \varphi \sigma_{\mu} \varphi
^{-1}$, and hence this composition is a homothety in $\mathbb{A}'_1$
with center at $p'_0$ that satisfies:

$$\varphi \sigma_{\mu} \varphi ^{-1}(p'_1)= \varphi
\sigma_{\mu}(p_1)= \varphi (p_{\mu})= p'_{h(\mu)} \ . $$

We conclude $h(\lambda \mu)= h(\lambda)h(\mu)$ by using the equality
$\varphi \sigma_{\mu} = \sigma '_{h(\mu)} \varphi$:
\begin{align*}
& \varphi \sigma_{\mu} (p_\lambda) = \varphi (p_{\lambda \mu})= p'_{h(\lambda \mu)} \\
& \sigma '_{h(\mu)} \varphi (p_{\lambda})= \sigma ' _{h(\mu)}
(p'_{h(\lambda)}) = p'_{h(\lambda)h(\mu)} \ .
\end{align*}

\hfill $\square$

\subsubsection*{The outer automorphism of $S_6$}

As an application of Theorem \ref{fundamental}, let us construct an
automorphism of $S_6$ which is not inner; i.e., which does not
coincide with conjugation by an element of $S_6$.

Let $k$ be the field with 5 elements. The projective line over $k$
has 6 elements, so that there are $6\cdot 5\cdot 4$ homographies and
$\mathrm{PGl}_1$ is a subgroup of $S_6$ of index 6. Since the
identity is the unique automorphism of the field
$k=\mathbb{Z}/5\mathbb{Z}$, Theorem \ref{fundamental} states that
$\mathrm{PGl}_1$ coincides with its normalizer in $S_6$; hence it
has 6 conjugated subgroups $H_1=\mathrm{PGl}_1,H_2,\ldots
,H_6$.\label{s6}

Any permutation $\tau \in S_6$ of the projective line
defines, by conjugation, a permutation $F(\tau )$ of this set
$\{ H_1,\ldots ,H_6\} $. Thus, we obtain an automorphism
$$
S_6 = \Biy ( \PP_1) \xrightarrow{\quad F \quad} {\rm Perm} (\{ H_1,
\ldots , H_6 \} ) = S_6
$$ such that $F(\mathrm{PGl}_1)$ is contained in the stabilizer of $H_1$; i.e., in
the subgroup of all
permutations fixing the element $H_1$ (in fact they coincide, since
both subgroups have index 6).

As the image of a stabilizer under an inner automorphism is the
stabilizer of another point, we conclude that $\, F\,$ cannot be
inner, since no point of the projective line is fixed by the group
of all homographies $\, \PGl_1$.

\subsection{Affine Geometry}

Let $\AA_n$ be the complementary of an hyperplane $H$ on $\PP_n$,
and let $\mathrm{Aff}_n$ be the group of affinites; i.e., the group
of projectivities $\varphi \colon \PP_n \to \PP_n$ such that
$\varphi (H) = H$.

Linear subvarieties of $\AA_n$ (affine lines, affine planes,...) are
defined as the restriction of affine subvarieties on $\PP_n$.

\begin{lemma}\label{planoafin} Let $\varphi \colon \AA_n \to \AA_n$ be a collineation
of an affine space of dimension $n>1$ over a
field $k\neq \mathbb{F}_2$.

If $\Pi \subset \mathbb{A}_n$ is an affine plane, then
$\varphi(\Pi)$ is contained in some affine plane.
\end{lemma}

\demo Let $L_1,L_2$ be two lines in $\Pi$  intersecting at a point $z$. Their images
$\varphi(L_1)$ and $\varphi(L_2)$ are lines with one point in common, so both lie in some affine plane $\Pi '$.

For any point $ p \in \Pi {\rm -} ( L_1 \cup L_2)$, let $L_p \subset
\Pi$ be any line passing through $p$, not parallel neither to $L_1$
nor to $L_2$, and such that $z \notin L_p$. As $L_p$ intersects both
$L_1$ and $L_2$, it follows that $\varphi(p) \in \varphi
(L_p)\subset \Pi '$.

\hfill $\square$

\medskip

If the base field has 2 elements, then there not exist 3
\emph{different} affine collinear points; in that case all
bijections of $\AA _n$ are in fact collineations.

\begin{definition}
A Staudt projectivity $\Phi \colon \PP_n \to \PP_n$ is called a {\it
Staudt affinity} if it preserves the hyperplane at infinity; i.e.,
if $\Phi (H) = H$.
\end{definition}

\begin{lemma}\label{colineacionesafines}
Let $\varphi \colon \AA_n \to \AA_n$ be a collineation of an affine
space of dimension $n>1$ over a field $k\neq \mathbb{F}_2$.

There exists a unique Staudt affinity $\Phi \colon \PP_n \to \PP_n$
such that $\varphi = \Phi_{|\AA _n}$.
\end{lemma}

\demo We first define $\Phi$ on the points at infinity: if $p \in
H$, let $\Phi (p)$ be the point at infinity of $\varphi (L_p)$,
where $L_p$ is any line passing through $p$. If $L'_p$ is another
line parallel to $L_p$, Lemma \ref{planoafin} shows that $\varphi
(L_p)$ and $\varphi (L'_p)$ are contained in some affine plane;
since they do not intersect, they are parallel and both have the
same point at infinity $\Phi (p)$.

Then, it is enough to check that $\Phi$ is a collineation on the
hyperplane at infinity (for collineations correspond with Staudt
projectivities, \cite{Artin} Thm 2.26): if three points $p_1,
p_2,p_3$ of the infinity are collinear, they lie in the direction of
some affine plane $\Pi $. By Lemma \ref{planoafin}, $\varphi (\Pi )$
is contained in some affine plane $\Pi '$, so that $\Phi (p_1),\Phi
(p_2),\Phi (p_3) \in H$ are in the direction of $\Pi '$ and are thus
collinear.

\hfill $\square$

\medskip
Let us now characterize collinear points in terms of affinities.
Consider:

$$ \mathsf{A}_{p_1,p_2}(p_3) :=\{ p \in \mathbb{A}_n{\rm -}\{p_3\} \colon \
\exists \ \varphi \in \mathrm{Aff}_n \ \ \varphi (p_1)=p_1,\ \varphi
(p_2)=p_2,\ \varphi (p_3)=p \, \} \ . $$

Since any affinity fixing $p_1$ and $p_2$ also fixes the point at
infinity of the line passing through them, it has to be the identity
on such line. Hence,

\begin{lemma}\label{afinrecta}
Three different points $p_1,p_2,p_3 \in \mathbb{A} _n$ are collinear if and only if
$$\mathsf{A}_{p_1,p_2}(p_3)=
\mathsf{A}_{p_1,p_3}(p_2)= \mathsf{A}_{p_2,p_3}(p_1)=\emptyset \ .
$$
\end{lemma}

\begin{thm}\label{afin}
The group of automorphisms of the Affine Geometry $(\mathbb{A}_n,
\mathrm{Aff}_n)$ over a field $k\neq \mathbb{F}_2$ is the group of
Staudt affinities, for any $n \geq 1$.

\end{thm}

\demo On the one hand, it is trivial to check that Staudt affinities
are indeed automorphisms, arguing as in the proof of Theorem
\ref{fundamental}.

On the other hand, let us prove that any automorphism of
$(\mathbb{A}_n, \mathrm{Aff}_n)$ is indeed a Staudt affinity:

\medskip
$n\geq 2.$ If $\varphi \in \mbox{Biy}(\AA_n)$ is in the normalizer
of $\mathrm{Aff}_n$, it preserves affinities and so $\varphi
(\mathsf{A}_{p_1,p_2}(p_3)) = \mathsf{A}_{\varphi ( p_1),\varphi
(p_2)}(\varphi (p_3))$. By Lemma \ref{afinrecta}, it follows that
$\varphi $ is a collineation of $\AA_n$, and hence it coincides with
the restriction to $\AA _n $ of a unique Staudt affinity  (Lemma
\ref{colineacionesafines}).

\medskip
$n=1$. Let $\varphi \in \mbox{Biy}(\AA _1)$ be in the normalizer of
$\mathrm{Aff}_1$, and fix a projective reference $( p_0,p_\infty,p_1
)$ such that $\AA _1=\PP _1 {\rm -} \{ p_\infty \} $. Via the
induced bijection $\AA _1\simeq k$, let $p_{\lambda }$ denote the
point corresponding to $\lambda \in k$. Consider $p'_0=\varphi
(p_0)$ and $p'_1=\varphi (p_1)$ as another origin and unit point,
and let $p'_{\lambda }$ denote the point with coordinate $\lambda
\in k$ via the corresponding bijection $\AA _1\simeq k$.

The composition $k \simeq \mathbb{A}_1 \xrightarrow{\varphi}
\mathbb{A}_1 \simeq k $ defines a bijection $h \colon k\to k$ such
that
$$\varphi(p_\lambda)=p'_{ h(\lambda)} \ .$$

A similar proof to that of Lemma \ref{automorfismo} shows that $h$
is an automorphism of the field $k$.

Then, it is easy to check that the bijection $\Phi \colon \PP_1 \to
\PP_1$ defined as $\Phi |_{\AA _1} :=\varphi$, $\Phi (p_\infty )
:=p_\infty$, is a Staudt afinity (with associated automorphism $h$).

\hfill $\square$

\subsection{Euclidean Geometry}

Let $\AA_n (\RR )$ be the affine space over the real numbers,
endowed with a positive definite, non-singular metric $g$ on its
vector space of directions. Let $\mbox{Mot}_n$ denote the group of
motions; i.e., those affinities $\varphi \colon \AA_n (\RR) \to
\AA_n (\RR)$ whose tangent linear map $\vec{\varphi}$ preserves the
metric $g$.

\begin{lemma}\label{perp} If an affinity $\varphi \colon \AA_n (\RR) \to
\AA_n (\RR)$ preserves motions (i.e., $\varphi \Mot_n \varphi^{-1} =
\Mot_n$), then it maps perpendicular lines into perpendicular lines.
\end{lemma}

\demo Assume there exists perpendicular lines $L_1,L_2$ such that
$\varphi(L_1)$ and $ \varphi(L_2)$ are not perpendicular.

Let $H_1$ be the hyperplane perpendicular to $L_2$ passing through
$L_1$. The symmetry $\sigma$ with respect to $H_1$ is a motion such
that $\varphi \sigma \varphi^{-1}$ is not a motion, for this
composition (which is not the identity), fixes the hyperplane
$\varphi(H_1)$ and preserves the oblique line $\varphi(L_2)$.

\hfill $\square$

\begin{definition}
An affinity $\varphi \colon \AA_n ( \RR) \to \AA_n ( \RR)$ is a {\it
similarity} if there exists $\lambda \in \RR$ such that
$\vec{\varphi}_* g = \lambda^2 g$.
\end{definition}

Observe that motions are a particular case of similarities, for
which $\lambda = 1$.

The Fundamental Theorem of Euclidean Geometry (\cite{Artin})
characterizes similarities as those affinities $\varphi \colon \AA_n
(\RR) \to \AA_n (\RR)$ that map perpendicular lines into
perpendicular lines.

\medskip

Analogously to what is made in previous sections, collinear points
may be characterized in terms of motions, considering the sets:

$$ \mathsf{M}_{p_1,p_2}(p_3) :=\{p\in \mathbb{A} _n - \{ p_3 \} \colon \,
\exists \, \varphi \in \mbox{Mot}_n \ , \ \varphi (p_1)=p_1,\
\varphi (p_2)=p_2,\ \varphi (p_3)=p \   \} \ . $$

\begin{lemma}\label{euclideo}
Three different points $p_1,p_2,p_3 \in \mathbb{A}_n(\mathbb{R})$
are collinear if and only if
$$\mathsf{M}_{p_1,p_2}(p_3)=
\mathsf{M}_{p_1,p_3}(p_2)=\mathsf{M}_{p_2,p_3}(p_1)=\emptyset$$
\end{lemma}

\begin{thm}\label{Euclidean}
The group of automorphisms of the Euclidean Geometry $(\AA _n(\mathbb{R}),
\mathrm{Mot}_n)$ is the group of similarities, for any $n>1$.

\end{thm}

\demo If $\varphi \colon \AA_n (\RR) \to \AA_n (\RR)$ is a
similarity, it is routine to check that, for any motion $\tau \in
\mathrm{Mot}_n$, the composition $\varphi \tau \varphi^{-1}$ is also
a motion.

On the other hand, if $\varphi \in \mbox{Biy}(\AA_n (\mathbb{R}))$
is a bijection that preserves motions, then $\varphi
(\mathsf{M}_{p_1,p_2}(p_3)) = \mathsf{M}_{\varphi ( p_1),\varphi
(p_2)}(\varphi (p_3))$ and, using Lemma \ref{euclideo}, it follows
that $\varphi $ is a collineation of $\AA_n(\mathbb{R})$.

As the identity is the unique automorphism of the field
$\mathbb{R}$, Lemma \ref{colineacionesafines} then assures that
$\varphi$ is an affinity.

Finally, as $\varphi $ is an affinity that preserves motions, it
also preserves perpendicular lines (Lemma \ref{perp}) and we can
assure that $\varphi $ is indeed a similarity.

\hfill $\square$

\bigskip

The above theorem is false when $n=1$: a counterexample is any
bijection $\varphi \colon \mathbb{R} \rightarrow \mathbb{R}$ that
respects the addition, v. gr., any $\mathbb{Q}$-linear automorphism
$\varphi$ of $\mathbb{R}$: if $\phi \in \mbox{Mot}_1$, then $\phi
(x)= \pm x + b$ and it follows that $(\varphi ^{-1} \phi \varphi
)(x) = \varphi ^{-1}( \pm \varphi (x) + b) = \pm x + \varphi ^{-1}
(b)$, so that $\varphi$ belongs to the normalizer of
$\mathrm{Mot}_1$.

However, it is not difficult to show that any \textit{continuous}
automorphism $\varphi $ of the euclidean line $(\AA _1(\mathbb{R}),
\mathrm{Mot}_1)$ is a similarity: $\varphi (x)=ax+b$; that is, the
group of similarities is the normalizer of $\mathrm{Mot}_1$ in the
group of all homeomorphisms of $\AA _1(\mathbb{R})$.

\subsection*{Acknowledgments}

The authors thank Juan A. Navarro and Juan B. Sancho for generous
advice and for their explanations on the outer automorphism of
$S_6$.

The first author has been partially supported by the Ministerio de
Econom\'{\i}a y Competitividad of Spain and the European Science
Foundation; the second author has been partially supported by Junta
de Extremadura and FEDER funds.

\end{document}